# Skew-product representations of multidimensional Dunkl Markov processes

## Oleksandr Chybiryakov

*Laboratoire de Probabilités et Modèles Aléatoires, Université Pierre et Marie Curie, 175 rue du Chevaleret, F-75013 Paris.*
*E-mail: oleksandr.chybiryakov@m4x.org*

**Abstract.** In this paper we obtain skew-product representations of the multidimensional Dunkl processes which generalize the skew-product decomposition in dimension 1 obtained in L. Gallardo and M. Yor. Some remarkable properties of the Dunkl martingales. *Séminaire de Probabilités XXXIX*, 2006. We also study the radial part of the Dunkl process, i.e. the projection of the Dunkl process on a Weyl chamber.

**Résumé.** Dans cet article nous obtenons des produits semi-directs des processus de Dunkl multidimensionnels qui généralisent ceux obtenus en dimension 1 dans L. Gallardo and M. Yor. Some remarkable properties of the Dunkl martingales. In *Séminaire de Probabilités XXXIX*, 2006. Nous étudions les processus de Dunkl radiaux qui sont les projections des processus de Dunkl sur une chambre de Weyl.

## 1. Introduction

The study of the multidimensional Dunkl processes was originated in [18] and [20]. They were studied further in [7, 8, 9, 10]. These processes share some important properties with Brownian motion: for example, they are martingales and enjoy the chaotic representation property as well as the time-inversion property. To a Dunkl process we can associate two processes: its Euclidean norm and its radial part. The Euclidean norm of the Dunkl process is in fact a Bessel process and its radial part is a continuous Markov process taking values in a Weyl chamber $C$. Note that a particular case of the radial part process – Brownian motion in a Weyl chamber – is studied in [2].

It is well known that Brownian motion (and more generally any rotational invariant diffusion) can be decomposed as the skew-product of a certain radial process and a time-changed spherical Brownian motion (see [6, 17]). In general, the following question can be raised: suppose given a group $G$ acting (not necessarily transitively) on $\mathbb{R}^n$ and a Markov process on $\mathbb{R}^n$, "invariant" by the action of $G$. Does this lead to a certain skew-product decomposition of the Markov process? In this paper we answer this question in the case of the Dunkl processes and the group $W$ – the group of symmetries of $\mathbb{R}^n$. We also study $X^W$ – the radial part of the Dunkl process noting some analogies between $X^W$ and Bessel processes.

The paper is organized as follows. In Section 2 we introduce some notations about martingale problems and extended infinitesimal generators for Markov processes. In Section 3 we give the definition of the Dunkl process and its radial part. We also present the invariance property of the Dunkl process under the action of the group $O(\mathbb{R}^n)$.

Section 4 contains the study of the radial part of the Dunkl process. We give the condition on the parameters of $X^W$ which guarantees that it does not hit the walls of a Weyl chamber $C$. This is done by applying a famous argument due to McKean (see Problem 7, p. 47 in [16]). As a part of the proof we find harmonic functions for $X^W$ which are analogs of harmonic functions for Bessel processes, i.e. the adequate power or logarithm function. Under the condition that $X^W$ does not hit the walls of $C$ it can be characterized as the unique solution of a stochastic integral equation or as the unique solution of the corresponding martingale problem (if $X^W$ hits the walls of $C$ it can still be characterized as a unique solution of the stochastic integral equation up to $T_0$ – the first time it hits the walls of $C$).

In Section 5 we generalize the skew-product decomposition given in [10] to multidimensional Dunkl processes. From [9] it is known that the Dunkl process can jump only in the directions given by the roots of an associated root system. In Theorem 19 we construct the Dunkl process from its radial part by adding the jumps in these directions one by one. In order to add the jumps we do the time-change then add jumps in the fixed direction and then do the inverse time-change. In order to add the jumps in a fixed direction we use the perturbation of generators technique from [5]. We also see that under a certain invariance condition for the extended infinitesimal generator this perturbation is given by a skew-product with a Poisson process. From Lemma 18 and Theorem 19 we obtain that the Dunkl process can be characterized as the unique solution of a martingale problem. Note a similar result from [20], which is recalled in Theorem 4. Our result holds for a more general class of processes – the two-parameter Dunkl processes (see [10] in dimension 1, [7] and [15]), but we need to impose the condition which guarantees that the radial part of the Dunkl process does not hit the wall of $C$.

## 2. Notations

We will constantly use martingale problems. Therefore, we need to recall here some definitions from [5] which will be used later.

Let $(S,d)$ be a metric space and $\mathbf{D}_S[0,\infty)$ ($\mathbf{C}_S[0,\infty)$) the space of right continuous (continuous) functions from $[0,\infty)$ into $(S,d)$ having left limits. $\mathcal{P}(S)$ denotes the set of Borel probability measures on $S$. For any $x \in S$, $\delta_x$ denotes the element of $\mathcal{P}(S)$ with unit mass at $x$. Let $L$ be the space of all measurable functions on $S$. $\mathcal{A}$ is a linear mapping whose domain $\mathcal{D}(\mathcal{A})$ is a subspace of $L$ and whose range $\mathcal{R}(\mathcal{A})$ lies in $L$. Typically $\mathcal{D}(\mathcal{A})$ will be $\mathbf{C}_K^\infty(S)$ – the space of infinitely differentiable functions on $S$ with compact support.

Let $X$ be a measurable stochastic process with values in $S$ defined on some probability space $(\Omega, \mathcal{F}, P)$. For $\nu \in \mathcal{P}(S)$ we say that $X$ is a *solution of the* $\mathbf{D}_S[0,\infty)$ *martingale problem* $(\mathcal{A}, \nu)$ if $X$ is a process with sample paths in $\mathbf{D}_S[0,\infty)$, $\mathbb{P}(X_0 \in \cdot) = \nu(\cdot)$, and for any $u \in \mathcal{D}(\mathcal{A})$

$$u(X_t) - u(X_0) - \int_0^t \mathcal{A}u(X_s)\,ds$$

is an $(\mathcal{F}_t^X)$-martingale, where $\mathcal{F}_t^X := \mathcal{F}\{X_s, s \leq t\}$, where in this paper $\mathcal{F}\{X_s, s \leq t\}$ indicates the sigma-field generated by variables $X_s$, $s \leq t$. Let $U$ be an open set of $S$ and $X$ be a process with sample paths in $\mathbf{D}_S[0,\infty)$. Define the $(\mathcal{F}_t^X)$-stopping time

$$\tau := \inf\{t \geq 0 | X_t \notin U \text{ or } X_{t-} \notin U\}.$$

Then $X$ is a *solution of the stopped* $\mathbf{D}_S[0,\infty)$ *martingale problem* $(\mathcal{A}, \nu, U)$ if $P(X_0 \in \cdot) = \nu(\cdot)$, $X_\cdot = X_{\cdot \wedge \tau}$ a.s., and for any $u \in \mathcal{D}(\mathcal{A})$

$$u(X_t) - u(X_0) - \int_0^{t \wedge \tau} \mathcal{A}u(X_s)\,ds$$

is an $(\mathcal{F}_t^X)$-martingale. If there exists a unique solution of a (stopped) martingale problem we will say that the (stopped) martingale problem is well-posed.

In what follows it will be convenient to work with the extended infinitesimal generator (see [14, 21], VII). We recall here the definition from ([21], VII).

If $X$ is a Markov process with respect to $(\mathcal{F}_t)$, a Borel function $f$ is said to belong to the domain $\mathbb{D}_A$ of the *extended infinitesimal generator* (or *extended generator*) if there exists a Borel function $g$ such that, a.s., $\int_0^t |g(X_s)|\,\mathrm{d}s < +\infty$ for every $t$, and

$$f(X_t) - f(X_0) - \int_0^t g(X_s)\,\mathrm{d}s$$

is a $(\mathcal{F}_t, P_x)$-right-continuous martingale for every $x$.

## 3. Preliminaries

Let $x \cdot y$ denote the usual scalar product for $x$ and $y$ on $\mathbb{R}^n$. For any $\alpha \in \mathbb{R}^n \setminus \{0\}$, $\sigma_\alpha$ denotes the reflection with respect to the hyperplane $H_\alpha \subset \mathbb{R}^n$ orthogonal to $\alpha$. For any $x \in \mathbb{R}^n$, it is given by

$$\sigma_\alpha(x) = x - 2\frac{\alpha \cdot x}{\alpha \cdot \alpha}\alpha.$$

For our purposes we need the following definition (see, for example, [19]):

**Definition 1.** *Let $R \subset \mathbb{R}^n \setminus \{0\}$ be a finite set. Then $R$ is called a (restricted) root system, if*

1. $R \cap \mathbb{R}\alpha = \{\pm\alpha\}$ *for all* $\alpha \in R$;
2. $\sigma_\alpha(R) = R$ *for all* $\alpha \in R$.

*The subgroup $W \subset O(\mathbb{R}^n)$ which is generated by the reflections $\{\sigma_\alpha \mid \alpha \in R\}$ is called the Weyl reflection group associated with $R$.*

One can prove that for any root system $R$ in $\mathbb{R}^n$, the reflection group $W$ is finite and the set of reflections contained in $W$ is exactly $\{\sigma_\alpha, \alpha \in R\}$ (see [19]).

**Example 2 (Root system of type $A_{n-1}$).** Let $e_1, \ldots, e_n$ be the standard basis vectors of $\mathbb{R}^n$, then

$$R = \{\pm(e_i - e_j), 1 \le i < j \le n\}$$

*is a root system in $\mathbb{R}^n$.*

**Example 3 (Root system of type $B_n$).** $R = \{\pm e_i, 1 \le i \le n, \pm(e_i \pm e_j), 1 \le i < j \le n\}$ *is a root system on $\mathbb{R}^n$.*

Without loss of generality we will suppose that if $R$ is a root system in $\mathbb{R}^n$ then for any $\alpha \in R$, $\alpha \cdot \alpha = 2$, so that for any $x \in \mathbb{R}^n$

$$\sigma_\alpha(x) = x - (\alpha \cdot x)\alpha. \tag{1}$$

For any given root system $R$ take $\beta \in \mathbb{R}^n \setminus \bigcup_{\alpha \in R} H_\alpha$, then $R_+ = \{\alpha \in R \mid \alpha \cdot \beta > 0\}$ is its positive subsystem. For any $\alpha \in R$ either $\alpha \in R_+$ or $-\alpha \in R_+$. Of course the choice of $R_+$ is not unique.

From [8, 20] the *Dunkl processes* $X^{(k)}$ are a family of càdlàg Markov processes with extended generators $\mathcal{L}_k$ where, for any $u \in \mathbf{C}_K^2(\mathbb{R}^n)$, $\mathcal{L}_k u$ is given by

$$\mathcal{L}_k u(x) = \frac{1}{2}\Delta u(x) + \sum_{\alpha \in R_+} k(\alpha)\left[\frac{\nabla u(x) \cdot \alpha}{x \cdot \alpha} - \frac{u(x) - u(\sigma_\alpha x)}{(x \cdot \alpha)^2}\right], \tag{2}$$

596 O. Chybiryakovwhere $k$ is a nonnegative multiplicity function invariant by the finite reflection group $W$ associated with $R$, i.e. $k: R \to [0, +\infty)$ and $k \circ w \equiv k$, for any $w \in W$. It is simple to see that $\mathcal{L}_k$ does not depend on the choice of $R_+$. In what follows suppose that $X_0 \in \mathbb{R}^n \setminus \bigcup_{\alpha \in R} H_\alpha$ a.s.

The semi-group densities of the Dunkl process with the generator (2) are given by

$$p_t^{(k)}(x, y) = \frac{1}{c_k t^{\gamma + n/2}} \exp\left(-\frac{|x|^2 + |y|^2}{2t}\right) D_k\left(\frac{x}{\sqrt{t}}, \frac{y}{\sqrt{t}}\right) \omega_k(y), \quad x, y \in \mathbb{R}^n, \tag{3}$$

where $D_k(u, v) > 0$ is the Dunkl kernel (for the properties of the Dunkl kernel see [19]), $\gamma := \sum_{\alpha \in R_+} k(\alpha)$, $\omega_k(y) = \prod_{\alpha \in R_+} |\alpha \cdot y|^{2k(\alpha)}$ and $c_k = \int_{\mathbb{R}^n} e^{-|x|^2/2} \omega_k(x) \, dx$ (see [18]).

We will need the following result from [20].

**Theorem 4.** *Let $(P_t)_{t \geq 0}$ be the semigroup of the Dunkl process given by its density (3) and let $\mathcal{L}_k$ be given by (2). Let $(X_t)_{t \geq 0}$ be a càdlàg process on $\mathbb{R}^n$ such that its Euclidean norm process $(\|X_t\|)_{t \geq 0}$ on $[0, +\infty[$ is continuous. Then the following statements are equivalent.*

1. *$(X_t)_{t \geq 0}$ is the Dunkl process associated with the semigroup $(P_t)_{t \geq 0}$.*
2. *For any $f \in \mathbf{C}^2(\mathbb{R}^n)$*

$$(M_t^f)_{t \geq 0} := \left(f(X_t) - f(X_0) - \int_0^t \mathcal{L}_k f(X_s) \, ds\right)_{t \geq 0}$$

*is a local martingale.*
3. *For any $f \in \mathbf{C}_K^2(\mathbb{R}^n)$, $(M_t^f)_{t \geq 0}$ is a martingale.*

Consider now a fixed Weyl chamber $C$ of the root-system $R$ which is a connected component of $\mathbb{R}^n \setminus \bigcup_{\alpha \in R} H_\alpha$. The space $\mathbb{R}^n / W$ of $W$-orbits in $\mathbb{R}^n$ can be identified to $\overline{C}$, i.e. there exists a homeomorphism $\phi: \mathbb{R}^n / W \to \overline{C}$. Denote $X_t^W := \pi(X_t^{(k)})$ – a *radial part of the Dunkl process* (or a *radial Dunkl process*), where

$$\pi := \phi \circ \pi_1 \tag{4}$$

and $\pi_1: \mathbb{R}^n \to \mathbb{R}^n / W$ denotes the canonical projection. From [8], $X^W$ is a Markov process with extended generator

$$\mathcal{L}_k^W u(x) = \frac{1}{2} \Delta u(x) + \sum_{\alpha \in R_+} k(\alpha) \frac{\nabla u(x) \cdot \alpha}{x \cdot \alpha}, \tag{5}$$

for any $u \in \mathbf{C}_0^2(\overline{C})$, such that $\nabla u(x) \cdot \alpha = 0$ for $x \in H_\alpha$, $\alpha \in R_+$. The semigroup densities of $X^W$ are of the form

$$p_t^W(x, y) = \frac{1}{c_k t^{\gamma + n/2}} \exp\left(-\frac{|x|^2 + |y|^2}{2t}\right) D_k^W\left(\frac{x}{\sqrt{t}}, \frac{y}{\sqrt{t}}\right) \omega_k(y), \quad x, y \in \overline{C}, \tag{6}$$

where

$$D_k^W(u, v) = \sum_{w \in W} D_k(u, wv). \tag{7}$$

From [20] the Dunkl process $(X_t)_{t \geq 0}$ is a Feller process and its Euclidean norm $(\|X_t\|)_{t \geq 0}$ is a Bessel process of index $\gamma + n/2 - 1$. It is easy to see that $(X_t^W)_{t \geq 0}$ defined above is also a Feller process.

In order to state the next result, we denote the trajectory of the Dunkl process starting at $x$ by $(X_t^{(k,R,x)})_{t \geq 0}$. The following proposition is an analog for the Dunkl–Markov processes of the rotational invariance property of the Brownian motion in $\mathbb{R}^n$.

**Proposition 5.** *Let $(X_t^{(k,R,x)})_{t\geq 0, x\in\mathbb{R}^n}$ be the Dunkl process, $O(\mathbb{R}^n)$ – the orthogonal group of $\mathbb{R}^n$. Then for any $\theta \in O(\mathbb{R}^n)$*

$$(\theta X_t^{(k,R,x)})_{t\geq 0, x\in\mathbb{R}^n} \stackrel{(d)}{=} (X_t^{(k_\theta, \theta R, \theta x)})_{t\geq 0, x\in\mathbb{R}^n},$$

*where $k_\theta : \theta R \ni \alpha \to k({}^t\theta \alpha)$, ${}^t\theta$ is the transpose of $\theta$, $\theta R = \{\theta \alpha, \alpha \in R\}$.*

**Proof.** It will be convenient to denote $\mathcal{L}_k^R$ the extended generator of $(X_t^{(k,R,x)})_{t\geq 0, x\in\mathbb{R}^n}$ given by (2). By Theorem 4 it is enough to prove that for any $u \in \mathbf{C}^2(\mathbb{R}^n)$ and $x \in \mathbb{R}^n$

$$\mathcal{L}_{k_\theta}^{\theta R} u(\theta x) = \mathcal{L}_k^R (u \circ \theta)(x).$$

From (1) one has

$$(u \circ \theta)(\sigma_\alpha(x)) = u(\theta x - (\alpha \cdot x)\theta\alpha) = u(\theta x - (\theta\alpha \cdot \theta x)\theta\alpha) = u(\sigma_{\theta\alpha}(\theta x))$$

and

$$\nabla(u \circ \theta)(x) \cdot \alpha = ({}^t\theta \nabla u(\theta x)) \cdot \alpha = \nabla u(\theta x) \cdot \theta\alpha.$$

Then

$$\mathcal{L}_k^R(u \circ \theta)(x) = \frac{1}{2}\Delta u(\theta x) + \sum_{\theta\alpha \in \theta R_+} k(\alpha) \left[ \frac{\nabla u(\theta x) \cdot \theta\alpha}{\theta x \cdot \theta\alpha} - \frac{u(\theta x) - u(\sigma_{\theta\alpha}(\theta x))}{(\theta x \cdot \theta\alpha)^2} \right]$$

$$= \mathcal{L}_{k_\theta}^{\theta R} u(\theta x). \qquad \square$$

## 4. Study of the Markov process $X^W$

Consider the Weyl chamber $C = \{x \in \mathbb{R}^n \mid x \cdot \alpha > 0, \alpha \in R_+\}$. We now study the radial part of the Dunkl process, i.e. the Markov process $(X_t^W)_{t\geq 0}$ with extended generator given by (5). $(X_t^W)_{t\geq 0}$ is a continuous Feller process with its values in $\overline{C}$.

Denote

$$\bar{\omega}_k(x) = \prod_{\alpha \in R_+} (\alpha \cdot x)^{k(\alpha)}, \qquad (8)$$

then for any $u \in \mathbf{C}_0^2(\overline{C})$, such that $\nabla u(x) \cdot \alpha = 0$ for $x \in H_\alpha$, $\alpha \in R_+$, and $x \in C$

$$\mathcal{L}_k^W u(x) = \frac{1}{2}\Delta u(x) + \sum_{\alpha \in R_+} k(\alpha) \frac{\nabla u(x) \cdot \alpha}{x \cdot \alpha}$$

$$= \frac{1}{2}\Delta u(x) + \nabla u(x) \cdot \nabla \log \bar{\omega}_k(x).$$

The following proposition gives a condition on the function $k: R \to \mathbb{R}_+$ in the definition of the extended generator (5), which ensures that the process $X^W$, such that $X_0^W \in C$ a.s., never touches $\partial C$.

**Proposition 6.** *Let $X^W$ be the radial Dunkl process, $X_0^W \in C$ a.s., with extended generator given by (5). Suppose that for any $\alpha \in R$, $k(\alpha) \geq \frac{1}{2}$. Define*

$$T_0 := \inf\{t > 0 \mid X_t^W \in \partial C\}, \qquad (9)$$

*then*

$$T_0 = +\infty \quad a.s.$$

**Remark 7.** *Suppose that for any $\alpha \in R, k(\alpha) \geq \frac{1}{2}$. From Proposition 6 $P_x(T_0 = +\infty) = 1$, where $P_x$ denotes the law of $X^W$ started at $x \in C$. Let $y \in \partial C$. For any $\varepsilon > 0$, using (6), one obtains that $P_y(X_\varepsilon^W \in H_\alpha \cap \overline{C}) = 0$. Since $\partial C \subset (\bigcup_{\alpha \in R} H_\alpha) \cap \overline{C}$*

$$P_y(X_\varepsilon^W \in C) = 1.$$

*Hence, by the same argument as in ([13], p. 162),*

$$P_y(X_t^W \in C, \forall \varepsilon < t < +\infty) = E_y(P_{X_\varepsilon^W}(X_t^W \in C, \forall 0 < t < +\infty)) = 1,$$

*where expectation $E_y$ is computed under the law $P_y$. Letting $\varepsilon \to 0$ one obtains*

$$P_y(T_0 = +\infty) = 1.$$

*Hence Proposition 6 is true if $X_0^W \in \partial C$.*

In order to prove this proposition we need the following two lemmas.

**Lemma 8.** *For $x \in C$ define*

$$\delta(x) := \prod_{\alpha \in R_+} (\alpha \cdot x)^{1-2k(\alpha)}, \tag{10}$$

*then $\delta$ is harmonic in $C$ for $\mathcal{L}_k^W$ given by (5), i.e. $\mathcal{L}_k^W \delta(x) = 0$ for any $x \in C$.*

**Remark 9.** *It is important in the following proof that the function $k$ in the definition of $\mathcal{L}_k^W$ is constant on the orbits of $R$ under the action of the associated Weyl group.*

**Proof of Lemma 8.** For $i = 1, \ldots, m$ let $R^i$ be the orbits of $R$ under the action of the associated Weyl group $W$. Denote $R_+^i = R^i \cap R_+$. Then for any $\{i_1, \ldots, i_l\} \subset \{1, \ldots, m\}$

$$\hat{R}_+ = R_+^{i_1} \cup \cdots \cup R_+^{i_l}$$

is also a positive subsystem of the root system $\hat{R} = R^{i_1} \cup \cdots \cup R^{i_l}$ in $\mathbb{R}^n$ (the two conditions in Definition 1 for $\hat{R}$ are easily verified). Without loss of generality it is enough to consider the irreducible case and in that case there are at most two orbits. Hence, it is enough to prove the lemma for $m = 2$. Denote

$$\pi_i(x) = \prod_{\alpha \in R_+^i} (\alpha \cdot x).$$

Since $k: R \to \mathbb{R}_+$ is constant on $R_+^i$ define

$$k_i := k(\alpha), \quad \alpha \in R_+^i,$$

then

$$\varpi_k(x) = \pi_1^{k_1}(x) \pi_2^{k_2}(x)$$

and

$$\delta(x) = \pi_1^{1-2k_1}(x) \pi_2^{1-2k_2}(x).$$

We want to prove that $\mathcal{L}_k^W \delta(x) = 0$ for any $x \in C$. From Theorem 4.2.6 in ([4], p. 140), for any root system $R$, if $\pi(x) = \prod_{\alpha \in R_+} (\alpha \cdot x)$, then

$$\Delta \pi(x) = 0.$$

Note that
$$\Delta \pi^{1-2k} + 2(\nabla \pi^{1-2k} \cdot \nabla \log \pi^k) = 0. \tag{11}$$

Indeed
$$\Delta \pi^{1-2k} = \nabla \cdot \nabla \pi^{1-2k} = \nabla \cdot ((1-2k)\pi^{-2k}\nabla \pi)$$
$$= (1-2k)[-2k\pi^{-2k-1}(\nabla \pi \cdot \nabla \pi) + \pi^{-2k}\Delta \pi]$$
$$= -2k(1-2k)\pi^{-2k-1}(\nabla \pi \cdot \nabla \pi).$$

On the other hand
$$2(\nabla \pi^{1-2k} \cdot \nabla \log \pi^k) = 2k(1-2k)\pi^{-2k-1}(\nabla \pi \cdot \nabla \pi)$$

and the result (11) follows.

Denote
$$\bar{\pi} := \pi_1^{1-2k_1}$$

and
$$\hat{\pi} := \pi_1^{k_1}.$$

We need to check that
$$A = \Delta(\pi_1^{1-2k_1}\pi_2^{1-2k_2}) + 2(\nabla(\pi_1^{1-2k_1}\pi_2^{1-2k_2}) \cdot \nabla \log(\pi_1^{k_1}\pi_2^{k_2})) = 0$$

or
$$\Delta(\bar{\pi}\pi_2^{1-2k_2}) + 2(\nabla(\bar{\pi}\pi_2^{1-2k_2}) \cdot \nabla \log(\hat{\pi}\pi_2^{k_2})) = 0.$$

One has
$$\Delta(\bar{\pi}\pi_2^{1-2k_2}) = \pi_2^{1-2k_2}\Delta\bar{\pi} + \bar{\pi}\Delta(\pi_2^{1-2k_2}) + 2(\nabla\bar{\pi} \cdot \nabla\pi_2^{1-2k_2})$$

and
$$\nabla(\bar{\pi}\pi_2^{1-2k_2}) \cdot \nabla \log(\hat{\pi}\pi_2^{k_2}) = (\bar{\pi}\nabla\pi_2^{1-2k_2} + \pi_2^{1-2k_2}\nabla\bar{\pi}) \cdot (\nabla \log \pi_2^{k_2} + \nabla \log \hat{\pi})$$
$$= (\bar{\pi}\nabla\pi_2^{1-2k_2} \cdot \nabla \log \pi_2^{k_2}) + (\bar{\pi}\nabla\pi_2^{1-2k_2} \cdot \nabla \log \hat{\pi})$$
$$+ ((\pi_2^{1-2k_2}\nabla\bar{\pi}) \cdot \nabla \log \pi_2^{k_2}) + ((\pi_2^{1-2k_2}\nabla\bar{\pi}) \cdot \nabla \log \hat{\pi}),$$

then
$$A = \pi_2^{1-2k_2}\Delta\bar{\pi} + 2\pi_2^{1-2k_2}(\nabla\bar{\pi} \cdot \nabla \log \hat{\pi}) + \bar{\pi}\Delta(\pi_2^{1-2k_2}) + 2\bar{\pi}(\nabla\pi_2^{1-2k_2} \cdot \nabla \log \pi_2^{k_2})$$
$$+ 2(\nabla\bar{\pi} \cdot \nabla\pi_2^{1-2k_2}) + 2((\bar{\pi}\nabla\pi_2^{1-2k_2}) \cdot \nabla \log \hat{\pi}) + 2((\pi_2^{1-2k_2}\nabla\bar{\pi}) \cdot \nabla \log \pi_2^{k_2}).$$

But
$$\Delta\bar{\pi} + 2(\nabla\bar{\pi} \cdot \nabla \log \hat{\pi}) = 0$$

and
$$\Delta(\pi_2^{1-2k_2}) + 2(\nabla\pi_2^{1-2k_2} \cdot \nabla \log \pi_2^{k_2}) = 0,$$

hence

$$A = 2(\nabla\bar{\pi} \cdot \nabla\pi_2^{1-2k_2}) + 2((\bar{\pi}\nabla\pi_2^{1-2k_2}) \cdot \nabla\log\hat{\pi}) + 2((\pi_2^{1-2k_2}\nabla\bar{\pi}) \cdot \nabla\log\pi_2^{k_2})$$
$$= 2\bar{\pi}\pi_2^{1-2k_2}[(\nabla\log\bar{\pi} \cdot \nabla\log\pi_2^{1-2k_2}) + (\nabla\log\pi_2^{1-2k_2} \cdot \nabla\log\hat{\pi}) + (\nabla\log\bar{\pi} \cdot \nabla\log\pi_2^{k_2})].$$

But one has

$$(\nabla\log\bar{\pi} \cdot \nabla\log\pi_2^{1-2k_2}) = (\nabla\log\pi_1^{1-2k_1} \cdot \nabla\log\pi_2^{1-2k_2}) = (1-2k_1)(1-2k_2)(\nabla\log\pi_1 \cdot \nabla\log\pi_2),$$
$$(\nabla\log\pi_2^{1-2k_2} \cdot \nabla\log\hat{\pi}) = k_1(1-2k_2)(\nabla\log\pi_1 \cdot \nabla\log\pi_2),$$
$$(\nabla\log\bar{\pi} \cdot \nabla\log\pi_2^{k_2}) = k_2(1-2k_1)(\nabla\log\pi_1 \cdot \nabla\log\pi_2)$$

and

$$(\nabla\log\pi_1 \cdot \nabla\log\pi_2) = \frac{1}{\pi_1\pi_2}(\nabla\pi_1 \cdot \nabla\pi_2)$$
$$= \frac{1}{2\pi_1\pi_2}(\Delta(\pi_1\pi_2) - \pi_1\Delta\pi_2 - \pi_2\Delta\pi_1) = 0.$$

Hence $A = 0$ and the lemma is proven. $\square$

In the same way one can prove the following lemma.

**Lemma 10.** *For $x \in C$ define*

$$\bar{\delta}(x) := \prod_{\alpha \in R_+, k(\alpha) \neq 1/2} (\alpha \cdot x)^{1-2k(\alpha)} \log \prod_{\alpha \in R_+, k(\alpha) = 1/2} (\alpha \cdot x), \tag{12}$$

*then $\bar{\delta}$ is harmonic in $C$ for $\mathcal{L}_k^W$ given by (5).*

**Proof.** Suppose that $k(\alpha) = \frac{1}{2}$ for any $\alpha \in R$ then for any $u \in \mathbf{C}^2(C)$, $x \in C$

$$\mathcal{L}_{1/2}^W u(x) = \frac{1}{2}(\Delta u(x) + (\nabla u(x) \cdot \nabla\log\pi(x))),$$

where $\pi(x) = \prod_{\alpha \in R_+}(\alpha \cdot x)$. Then

$$\Delta\log\pi(x) = -\frac{1}{\pi^2(x)}(\nabla\pi(x) \cdot \nabla\pi(x)) + \frac{1}{\pi(x)}\Delta\pi(x) = -\frac{1}{\pi^2(x)}(\nabla\pi(x) \cdot \nabla\pi(x))$$

and

$$\Delta\log\pi(x) + (\nabla\log\pi(x) \cdot \nabla\log\pi(x)) = 0.$$

Denote

$$\bar{\pi}_1 := \prod_{\alpha \in R_+, k(\alpha) \neq 1/2} (\alpha \cdot x)^{1-2k(\alpha)} = \prod_{1 \leq i \leq m, k_i \neq 1/2} \pi_i^{1-2k_i},$$

$$\hat{\pi}_1 := \prod_{\alpha \in R_+, k(\alpha) \neq 1/2} (\alpha \cdot x)^{k(\alpha)} = \prod_{1 \leq i \leq m, k_i \neq 1/2} \pi_i^{k_i},$$

$$\bar{\pi}_2 := \log \prod_{\alpha \in R_+, k(\alpha) = 1/2} (\alpha \cdot x),$$

$$\hat{\pi}_2 := \prod_{\alpha \in R_+, k(\alpha) = 1/2} (\alpha \cdot x)^{1/2},$$

then $\bar{\omega}_k = \hat{\pi}_1 \hat{\pi}_2$ and $\bar{\delta} = \bar{\pi}_1 \bar{\pi}_2$. In order to prove Lemma 10 one needs to check that

$$\Delta(\bar{\pi}_1 \bar{\pi}_2) + 2(\nabla(\bar{\pi}_1 \bar{\pi}_2) \cdot \nabla \log(\hat{\pi}_1 \hat{\pi}_2)) = 0,$$

but it is equal to

$$\bar{\pi}_2[\Delta \bar{\pi}_1 + 2(\nabla \bar{\pi}_1 \cdot \nabla \log(\hat{\pi}_1))] + \bar{\pi}_1[\Delta \bar{\pi}_2 + 2(\nabla \bar{\pi}_2 \cdot \nabla \log(\hat{\pi}_2))]$$
$$+ 2(\nabla \bar{\pi}_1 \cdot \nabla \bar{\pi}_2) + 2\bar{\pi}_2(\nabla \bar{\pi}_1 \cdot \nabla \log \hat{\pi}_2) + 2\bar{\pi}_1(\nabla \bar{\pi}_2 \cdot \nabla \log \hat{\pi}_1)$$

and in the same way as in the proof of Lemma 8 one can see that it is equal to zero. $\square$

Using Lemmas 8 and 10 one obtains the following result, which will lead to a presentation of $X^W$ (see Corollaries 13 and 14).

**Lemma 11.** *Let $\varpi_k$ be defined by (8) and $k(\alpha) \geq \frac{1}{2}$ for any $\alpha \in R$. Suppose that $X_0 \in C$ a.s., then there exists a unique solution $X$ of the stochastic integral equation*

$$X_t = X_0 + \beta_t + \int_0^t \nabla \log \varpi_k(X_s) \, ds, \tag{13}$$

*where $(\beta_t)_{t \geq 0}$ is a Brownian motion on $\mathbb{R}^n$. Furthermore $\mathbb{P}(\forall t > 0, X_t \in C) = 1$.*

**Proof.** We will follow the proof of the similar argument given in ([1] Lemma 3.2). Since the function $\nabla \log \varpi_k \in \mathbf{C}^{\infty}(C)$ Eq. (13) has a unique (strong) maximal solution in $C$, defined up to time $\zeta$, where $\zeta$ is an explosion time or the exit time from $C$.

*Case 1* ($k(\alpha) \neq \frac{1}{2}$ for any $\alpha \in R$). From Lemma 8 the function $\delta$ defined by (10) is harmonic and positive on $C$. By Ito's formula one deduces that $\{\delta(X_t), t < \zeta\}$ is a positive local martingale, thus it converges a.s. when $t \to \zeta$. But $\delta = +\infty$ on $\partial C$, therefore $\|X_t\| \to +\infty$ when $t \to \zeta$. On the other hand by Ito's formula for $t < \zeta$

$$\|X_t\|^2 = \|X_0\|^2 + 2\int_0^t (X_s \cdot d\beta_s) + 2\int_0^t \gamma \, ds + tn$$
$$= \|X_0\|^2 + 2\int_0^t \|X_s\| \, d\tilde{\beta}_s + t(n + 2\gamma),$$

where $\gamma = \sum_{\alpha \in R_+} k(\alpha)$ and $\tilde{\beta}_t = \int_0^t \|X_s\|^{-1}(X_s \cdot d\beta_s)$, $t < \zeta$ is a real-valued Brownian motion up to time $\zeta$. This shows that $\|X_t\|^2$ is the square of a $(n + 2\gamma)$-dimensional Bessel process up to time $\zeta$ (started at $\|X_0\|^2 > 0$ a.s.). Since $\|X_t\|^2 \to +\infty$, by standard results (see [21], XI.1) $\zeta = +\infty$ a.s. This implies that $T_0 = +\infty$ a.s., where $T_0$ is defined by (9).

*Case 2* (There exists $\alpha \in R$ such that $k(\alpha) = \frac{1}{2}$). From Lemma 10 and Ito's formula, $\{\bar{\delta}(X_t), t < \zeta\}$ is a continuous local martingale. Let $A_t := \langle \bar{\delta}(X), \bar{\delta}(X) \rangle_t$ for $t < \zeta$ and $\tau_t := \inf\{s \geq 0 \mid A_s = t\}$, then by Theorem 1.7 in ([21], Chapter V) $B_t = \bar{\delta}(X_{\tau_t}) - \bar{\delta}(X_0)$ will be a real-valued Brownian motion up to time $A_\zeta$. If $\zeta < +\infty$ we have already seen that $\|X_t\|$ cannot tend to $+\infty$. Hence, if $t \to \zeta$, either $\bar{\delta}(X_t)$ tends to $+\infty$ or to $-\infty$. Therefore, either $B_t$ tends to $+\infty$ or to $-\infty$, when $t \to A_\zeta$, but that is impossible for Brownian motion, either because $A_\zeta < +\infty$, or because $A_\zeta = +\infty$ and $B_t \geq 0$ ($B_t \leq 0$) infinitely often when $t \to +\infty$. $\square$

**Proof of Proposition 6.** Let

$$C_m := \left\{ x \in C \,\Big|\, x \cdot \alpha > \frac{1}{m}, \alpha \in R, \|x\| < m \right\}. \tag{14}$$

Take $g_m \in \mathbf{C}^\infty(\mathbb{R}^n)$ such that $g_m \equiv \log \varpi_k$ on $C_{m+1}$ and $g_m \equiv 0$ on $\mathbb{R}^n \setminus C_{m+2}$. For any $u \in \mathbf{C}_K^\infty(\mathbb{R}^n)$, define $\mathcal{L}_m$ by

$$\mathcal{L}_m u(x) = \frac{1}{2}\Delta u(x) + (\nabla u(x) \cdot \nabla g_m(x)).$$

Let $(X_t^W)_{t \geq 0}$ be the radial part of the Dunkl process and $\hat{\tau}_m := \inf\{s > 0 \mid X_s^W \notin C_m\}$. For any $f \in \mathbf{C}_K^\infty(\mathbb{R}^n)$ there exists $\tilde{f} \in \mathbf{C}_K^\infty(C)$ such that $\tilde{f} \equiv f$ on $C_{m+1}$. Then from (5)

$$\tilde{f}(X_{t \wedge \hat{\tau}_m}^W) - \tilde{f}(X_0^W) - \int_0^{t \wedge \hat{\tau}_m} \mathcal{L}_k^W \tilde{f}(X_s^W)\,\mathrm{d}s = f(X_{t \wedge \hat{\tau}_m}^W) - f(X_0^W) - \int_0^{t \wedge \hat{\tau}_m} \mathcal{L}_m f(X_s^W)\,\mathrm{d}s \tag{15}$$

is a martingale.

Let $X$ be the solution of (13) such that

$$\mathbb{P}(X_0 \in \cdot) = \mathbb{P}(X_0^W \in \cdot) = \bar{\nu}(\cdot) \tag{16}$$

for $\bar{\nu} \in \mathcal{P}(C)$. Let $\bar{\tau}_m := \inf\{s > 0 \mid X_s \notin C_m\}$. Then from Lemma 11 one deduces that

$$\bar{\tau}_m \to +\infty \quad \text{a.s.} \tag{17}$$

Furthermore from (13) and Ito's formula one deduces that

$$f(X_{t \wedge \bar{\tau}_m}) - f(X_0) - \int_0^{t \wedge \bar{\tau}_m} \mathcal{L}_m f(X_s)\,\mathrm{d}s \tag{18}$$

is a martingale for $f \in \mathbf{C}_K^\infty(\mathbb{R}^n)$.

By Theorem 3.3 in ([5], p. 379) for any $\nu \in \mathcal{P}(\mathbb{R}^n)$ the $\mathbf{D}_{\mathbb{R}^n}[0,\infty)$ martingale problem $(\mathcal{L}_m, \nu)$ is well posed. Then by Theorem 6.1 in ([5], p. 216) for each $\nu \in \mathcal{P}(\mathbb{R}^n)$ the stopped $\mathbf{D}_{\mathbb{R}^n}[0,\infty)$ martingale problem $(\mathcal{L}_m, \nu, C_m)$ is well-posed. But from (15), (18) and (16) $X_{\cdot \wedge \hat{\tau}_m}^W$ and $X_{\cdot \wedge \bar{\tau}_m}$ are solutions of $(\mathcal{L}_m, \bar{\nu}, C_m)$. Hence

$$(X_{t \wedge \hat{\tau}_m}^W)_{t \geq 0} \stackrel{(d)}{=} (X_{t \wedge \bar{\tau}_m})_{t \geq 0} \tag{19}$$

and $\mathbb{P}(\hat{\tau}_m < t) = \mathbb{P}(\bar{\tau}_m < t)$, for any $t > 0$. Hence from (17)

$$\hat{\tau}_m \to +\infty \quad \text{a.s.}$$

and passing to the limit as $m \to \infty$ in (19) one obtains that

$$(X_t^W)_{t \geq 0} \stackrel{(d)}{=} (X_t)_{t \geq 0}$$

and

$$T_0 = +\infty \quad \text{a.s.,}$$

where $T_0$ is defined by (9). $\square$

Since $X_t^W = \pi(X_t)$, where $X_t$ is the Dunkl process one obtains

**Corollary 12.** *Let $(X_t)_{t \geq 0}$ be the Dunkl process, such that $X_0 \in \mathbb{R}^n \setminus \bigcup_{\alpha \in R_+} H_\alpha$ a.s., with extended generator given by (2). Suppose that, for any $\alpha \in R$, $k(\alpha) \geq \frac{1}{2}$. Define*

$$T_0 := \inf\left\{t > 0 \,\Big|\, X_t \in \bigcup_{\alpha \in R_+} H_\alpha\right\},$$

*then*

$$T_0 = +\infty \quad a.s.$$

The proof of Proposition 6 leads to

**Corollary 13.** *Let* $(X_t^W)_{t\geq 0}$ *be the radial Dunkl process, such that* $X_0 \in C$ *a.s., with extended generator given by* (5). *Suppose that for any* $\alpha \in R$ $k(\alpha) \geq \frac{1}{2}$. *Then* $X^W$ *is the unique solution to the stochastic integral equation*

$$X_t = X_0 + \beta_t + \int_0^t \nabla \log \varpi_k(X_s)\, \mathrm{d}s,$$

*where* $(\beta_t)_{t\geq 0}$ *is a Brownian motion in* $\mathbb{R}^n$ *and* $\varpi_k$ *is given by* (8). *This solution is strong in the sense of ([12], IV).*

*Denote* $\nu(\cdot) := \mathbb{P}(X_0^W \in \cdot)$, *then* $X^W$ *is a unique solution to the* $\mathbf{C}_C[0,\infty)$ *martingale problem* $(\mathcal{L}_k^W, \nu)$, *where* $\mathcal{L}_k^W$ *is the restriction of* $\mathcal{L}_k^W$ *on* $\mathbf{C}_K^\infty(C)$.

Suppose now that there exists $\alpha \in R$ such that $k(\alpha) < \frac{1}{2}$ and consider the radial Dunkl process. The preceding results extend to this case if one works till $T_0$ – the first time the process hits the walls of the Weyl chamber. Loosely speaking, here, $T_0$ is considered as if it were an "explosion time." One gets the following corollary.

**Corollary 14.** *Let* $(X_t^W)_{t\geq 0}$ *be the radial Dunkl process, such that* $X_0 \in C$ *a.s., with extended generator given by* (5). *Suppose that there exists* $\alpha \in R$ *such that* $k(\alpha) < \frac{1}{2}$. *Let* $T_0$ *be defined by* (9). *Then* $(X_t^W, t < T_0)$ *is the unique solution to the stochastic integral equation*

$$X_t = X_0 + \beta_t + \int_0^t \nabla \log \varpi_k(X_s)\, \mathrm{d}s, \quad t < T_0,$$

*where* $(\beta_t)_{t\geq 0}$ *is a Brownian motion in* $\mathbb{R}^n$ *and* $\varpi_k$ *is given by* (8). *This solution is strong in the sense of ([12], IV).*

**Proof.** Is analogous to the proof of Lemma 11 and Proposition 6. □

*Remark 15.* One can show that the result of Corollary 14 is true for $t \geq 0$ and $X^W$ "reflects instantly" when it reaches $\partial C$ (see [3]).

## 5. Skew-product decomposition of jumps

Let $X$ be a Dunkl process, $\pi: \mathbb{R}^n \to \overline{C}$ be given by (4). Since $X$ is càdlàg and $\pi(X)$ is continuous it will be useful to introduce the space $\mathbf{D}_{\mathbb{R}^n}^\pi[0,\infty)$ defined by

$$\mathbf{D}_{\mathbb{R}^n}^\pi[0,\infty) := \{\omega \in \mathbf{D}_{\mathbb{R}^n}[0,\infty) \mid \pi \circ \omega \in \mathbf{C}_{\overline{C}}[0,\infty)\}. \tag{20}$$

Let $E := \mathbb{R}^n \setminus \bigcup_{\alpha \in R} H_\alpha$. Denote

$$\mathbf{D}_E^\pi[0,\infty) := \{\omega \in \mathbf{D}_E[0,\infty) \mid \pi \circ \omega \in \mathbf{C}_C[0,\infty)\}.$$

We will say that $X$ is a *solution of the* $\mathbf{D}_E^\pi[0,\infty)$ $(\mathbf{D}_{\mathbb{R}^n}^\pi[0,\infty))$ *martingale problem* $(\mathcal{A}, \nu)$ if $X$ is a process with sample paths in $\mathbf{D}_E^\pi[0,\infty)$ $(\mathbf{D}_{\mathbb{R}^n}^\pi[0,\infty))$ and $X$ is a solution of the $\mathbf{D}_E[0,\infty)$ $(\mathbf{D}_{\mathbb{R}^n}[0,\infty))$ martingale problem $(\mathcal{A}, \nu)$.

Let $\lambda > 0$ and $\alpha \in R_+$. Let $\mathcal{D}(\mathcal{A}) \subset \mathbf{C}_K^\infty(E)$ and $\mathcal{R}(\mathcal{A}) \subset \mathbf{C}_K^\infty(E)$. Suppose that for any $\nu \in \mathcal{P}(\mathbb{R}^n)$ there exists $X$ – a solution of the $\mathbf{D}_{\mathbb{R}^n}^\pi[0,\infty)$ martingale problem for $(\mathcal{A},\nu)$. We will need the following construction from ([5], p. 256). Let $\Omega = \prod_{k=1}^\infty (\mathbf{D}_{\mathbb{R}^n}[0,\infty) \times [0,\infty))$ and $(X_k, \Delta_k)$ denote the coordinate random variables. Define $\mathcal{G}_k = \mathcal{F}(X_l, \Delta_l\colon l \leq k)$ and $\mathcal{G}^k = \mathcal{F}(X_l, \Delta_l\colon l \geq k)$. Then there is a probability distribution on $\Omega$ such that for each $k$ $X_k$ is a solution of the $\mathbf{D}_{\mathbb{R}^n}^\pi[0,\infty)$ martingale problem for $\mathcal{A}$, $\Delta_k$ is independent of $\mathcal{F}(X_1, \ldots, X_k, \Delta_1, \ldots, \Delta_{k-1})$ and exponentially distributed with parameter $\lambda$, and for $A_1 \in \mathcal{G}_k$ and $A_2 \in \mathcal{G}^{k+1}$,

$$\mathbb{P}(A_1 \cap A_2) = \mathbb{E}(\mathbb{I}_{A_1} \mathbb{P}[A_2 | X_{k+1}(0) = \sigma_\alpha(X_k(\Delta_k))]) \tag{21}$$

and $\mathbb{P}(X_1(0) \in \cdot) = \nu(\cdot)$. Define $\tau_0 = 0$, $\tau_k = \sum_{i=1}^k \Delta_i$, and $N_t = k$ for $\tau_k \leq t < \tau_{k+1}$. Note that $N$ is a Poisson process with parameter $\lambda$. Define

$$Y(t) = X_{k+1}(t - \tau_k), \quad \tau_k \leq t < \tau_{k+1}, \tag{22}$$

and $\mathcal{F}_t := \mathcal{F}_t^Y \vee \mathcal{F}_t^N$.

**Notation 16.** *We will denote $Y$ in (22) by $X *_\alpha N$.*

**Lemma 17.** *Let $\mathcal{D}(\mathcal{A}) \subset \mathbf{C}_K^\infty(E)$ and $\mathcal{R}(\mathcal{A}) \subset \mathbf{C}_K^\infty(E)$. Suppose that for any $\nu \in \mathcal{P}(E)$ there exists a solution $X$ of the $\mathbf{D}_E^\pi[0,\infty)$ martingale problem $(\mathcal{A},\nu)$. Then for any $\nu \in \mathcal{P}(E)$ there exists a Poisson process $N$ with parameter $\lambda$ such that $Y = X *_\alpha N$, where $*_\alpha$ is defined by (22), is a solution of the $\mathbf{D}_E^\pi[0,\infty)$ martingale problem $(\mathcal{A}_{\lambda,\alpha}, \nu)$, where for any $u \in \mathcal{D}(\mathcal{A})$ and $x \in E$*

$$\mathcal{A}_{\lambda,\alpha} u(x) = \mathcal{A}u(x) + \lambda(u(\sigma_\alpha x) - u(x)).$$

*Furthermore if for any $u \in \mathcal{D}(\mathcal{A})$ and $x \in E$*

$$\mathcal{A}(u \circ \sigma_\alpha)(x) = \mathcal{A}u(\sigma_\alpha(x)), \tag{23}$$

*then there exists a solution of $(\mathcal{A}_{\lambda,\alpha}, \nu)$ given by*

$$(Y_t)_{t \geq 0} := (\sigma_\alpha^{N_t} X_t)_{t \geq 0},$$

*where $N$ is a Poisson process with parameter $\lambda$ independent of $X$.*

**Proof.** We follow the proof of Proposition 10.2 in ([5], p. 256). Since $E$ is not complete one shall consider the $\mathbf{D}_{\mathbb{R}^n}^\pi[0,\infty)$ martingale problem $(\mathcal{A}, \tilde{\nu})$. For any $\tilde{\nu} \in \mathcal{P}(\mathbb{R}^n)$ let $\tilde{X}_0$ be such that $\mathbb{P}(\tilde{X}_0 \in \cdot) = \tilde{\nu}(\cdot)$. If $\tilde{\nu}(E) = 0$, then $\tilde{X}_t := \tilde{X}_0$, for any $t \geq 0$, is a solution of the $\mathbf{D}_{\mathbb{R}^n}^\pi[0,\infty)$ martingale problem $(\mathcal{A}, \tilde{\nu})$. If $\tilde{\nu}(E) > 0$ and $(X_t)_{t \geq 0}$ is a solution of the $\mathbf{D}_E^\pi[0,\infty)$ martingale problem $(\mathcal{A}, \bar{\nu})$ with $\bar{\nu}(\cdot) = \tilde{\nu}(\cdot)/\tilde{\nu}(E)$, $\bar{\nu} \in \mathcal{P}(E)$, then $\tilde{X}_t := X_t \mathbb{I}_{\{\tilde{X}_0 \in E\}} + \tilde{X}_0 \mathbb{I}_{\{\tilde{X}_0 \notin E\}}$ is a solution of the $\mathbf{D}_{\mathbb{R}^n}^\pi[0,\infty)$ martingale problem $(\mathcal{A}, \tilde{\nu})$. Taking $\mu(x, \cdot) = \delta_{\sigma_\alpha(x)}(\cdot)$ and using Proposition 10.2 in ([5], p. 256) one obtains that $Y = \tilde{X} *_\alpha N$ is a solution of the $\mathbf{D}_{\mathbb{R}^n}^\pi[0,\infty)$ martingale problem $(\mathcal{A}_{\lambda,\alpha}, \tilde{\nu})$. Let $\tilde{\nu} = \nu \in \mathcal{P}(E)$, then $\tilde{X} \equiv X$ is a solution of the $\mathbf{D}_E^\pi[0,\infty)$ martingale problem $(\mathcal{A}, \nu)$. Furthermore, from (21)

$$\mathbb{P}\left(\bigcap_{k \geq 1} \{X_{k+1}(0) = \sigma_\alpha(X_k(\Delta_k))\}\right) = 1$$

and

$$\mathbb{P}\left(\bigcap_{k \geq 1} \{Y(\tau_k) = \sigma_\alpha(Y(\tau_k-))\}\right) = 1.$$

Hence
$$\mathbb{P}\left(\bigcap_{k\geq 1}\{\pi(Y(\tau_k))=\pi(Y(\tau_k-))\}\right)=1.$$

Since $X_k$ has paths in $\mathbf{D}_E^\pi[0,\infty)$ for any $k\geq 1$ $\pi(Y)$ is a.s. continuous and $Y$ is a process with sample paths in $\mathbf{D}_E[0,\infty)$. Hence $Y$ is a solution of the $\mathbf{D}_E^\pi[0,\infty)$ martingale problem $(\mathcal{A}_{\lambda,\alpha},\nu)$.

Suppose now that (23) is true. Let $N$ be a Poisson process independent of $X$ and $\tau_k := \inf\{s \geq 0 \,|\, N_s = k\}$. Let $Y_t = \sigma_\alpha^{N_t} X_t$. Note that for any $u \in \mathbf{C}_K^\infty(E)$

$$M_t^u := u(X_t) - u(X_0) - \int_0^t \mathcal{A}u(X_s)\,\mathrm{d}s$$

is a $(\mathcal{F}_t^X)$-martingale.

Since $(\tau_i)_{i\geq 0}$ are independent from $M^u$, for any $k\geq 0$

$$u(X((t\vee \tau_{2k})\wedge \tau_{2k+1})) - u(X(\tau_{2k})) - \int_{\tau_{2k}}^{(t\vee \tau_{2k})\wedge \tau_{2k+1}} \mathcal{A}u(X_s)\,\mathrm{d}s \tag{24}$$

is a $(\mathcal{F}_t)$-martingale. From (23)

$$u(\sigma_\alpha X((t\vee \tau_{2k+1})\wedge \tau_{2k+2})) - u(\sigma_\alpha X(\tau_{2k+1})) - \int_{\tau_{2k+1}}^{(t\vee \tau_{2k+1})\wedge \tau_{2k+2}} \mathcal{A}u(\sigma_\alpha X_s)\,\mathrm{d}s \tag{25}$$

is a $(\mathcal{F}_t)$-martingale. Summing (24) and (25) over $k$ one gets that

$$u(\sigma_\alpha^{N_t} X_t) - u(X_0) - \int_0^t \mathcal{A}u(\sigma_\alpha^{N_s} X_s)\,\mathrm{d}s - \sum_{k=1}^{N_t}(u(\sigma_\alpha^{k+1} X_{\tau_k}) - u(\sigma_\alpha^k X_{\tau_k})) \tag{26}$$

is a $(\mathcal{F}_t)$-martingale. But

$$\int_0^t (u(\sigma_\alpha Y_{s-}) - u(Y_{s-}))\,\mathrm{d}(N_s - \lambda s) \tag{27}$$

also is a $(\mathcal{F}_t)$-martingale. Note that since $N$ is independent of $X$

$$\int_0^t (u(\sigma_\alpha Y_{s-}) - u(Y_{s-}))\,\mathrm{d}N_s = \sum_{k=1}^{N_t}(u(\sigma_\alpha Y_{\tau_k}) - u(Y_{\tau_k})).$$

Adding (26) and (27) and noting that: $\sigma_\alpha^k X_{\tau_k} = Y_{\tau_k}$, one gets that

$$u(Y_t) - u(X_0) - \left(\int_0^t \mathcal{A}u(Y_s)\,\mathrm{d}s + \lambda \int_0^t (u(\sigma_\alpha Y_s) - u(Y_s))\,\mathrm{d}s\right)$$

is a $(\mathcal{F}_t)$-martingale. □

**Lemma 18.** *Suppose that $k: R \to [\frac{1}{2}, +\infty)$ is a $W$-invariant multiplicity function, $l: R \to [0, +\infty)$, and for any $u \in \mathbf{C}_K^\infty(E)$, $x \in E$*

$$\mathcal{L}_{k,l} u(x) = \frac{1}{2}\Delta u(x) + \sum_{\alpha \in R_+} k(\alpha)\frac{\nabla u(x)\cdot \alpha}{x\cdot \alpha} + \sum_{\alpha \in R_+} l(\alpha)\frac{u(\sigma_\alpha x) - u(x)}{(x\cdot \alpha)^2}.$$

*Let $\pi: E \to \overline{C}$ be defined by (4). For any $\nu \in \mathcal{P}(E)$, if there exists a solution of the $\mathbf{D}_E^\pi[0,\infty)$ martingale problem $(\mathcal{L}_{k,l},\nu)$, then it is unique.*

**Proof.** Let $X$ be a solution of the $\mathbf{D}_E^\tau[0,\infty)$ martingale problem $(\mathcal{L}_{k,l}, \nu)$. Let us fix a Weyl chamber $C$ and a projection $\pi : E \to \overline{C}$, defined by (4). Denote $Y_t := \pi(X_t)$, then taking functions $u \in \mathbf{C}_K^\infty(E)$ such that $u(\sigma_\alpha x) = u(x)$ for any $\alpha \in R$ and $x \in E$ one sees that $u(X_\cdot) = u(Y_\cdot)$ and

$$u(Y_t) - u(Y_0) - \int_0^t \mathcal{L}_k^W u(Y_s)\, ds$$

is a martingale. Hence, from the results in the previous section $Y$ is a unique solution of the $\mathbf{C}_C[0,\infty)$ martingale problem $(\mathcal{L}_k^W, \hat{\nu})$, where $\hat{\nu}(\cdot) = \nu(\pi^{-1}(\cdot))$, i.e. $Y$ is a radial Dunkl process. Define

$$B_m := \bigcup_{w \in W} w(C_m),$$

where $C_m$ is defined by (14) and $T_m := \inf\{s \geq 0 | X_s \notin B_m \text{ or } X_{s-} \notin B_m\}$. Since $X_s \notin B_m \iff Y_s \notin C_m$ and $Y$ is a.s. continuous

$$T_m = \inf\{s \geq 0 \,|\, Y_s \notin C_m\} \to \infty \quad \text{as } m \to \infty \quad \text{a.s.} \tag{28}$$

and for any $t \geq 0$ and $m$ sufficiently large, such that $X_0 \in \overline{C}_m$,

$$X_{t \wedge T_m} \in \overline{C}_m. \tag{29}$$

As in the proof of Proposition 6 there exist bounded functions $g_m \in \mathbf{C}^\infty(\mathbb{R}^n)$, $f_{\alpha,m} \in \mathbf{C}^\infty(\mathbb{R}^n)$, such that $g_m \equiv \log \varpi_k$ on $B_{m+1}$ and $g_m \equiv 0$ on $\mathbb{R}^n \setminus B_{m+2}$, where $\varpi_k$ is defined by (8), and for any $\alpha \in R$ $f_{\alpha,m}(x) = \frac{1}{(x \cdot \alpha)^2}$ for $x \in B_{m+1}$ and $f_{\alpha,m} \equiv 0$ on $\mathbb{R}^n \setminus B_{m+2}$. For any $u \in \mathbf{C}_K^\infty(\mathbb{R}^n)$, define $\mathcal{A}_m$ by

$$\mathcal{A}_m u(x) = \frac{1}{2}\Delta u(x) + (\nabla u(x) \cdot \nabla g_m(x)) + \sum_{\alpha \in R_+} l(\alpha) f_{\alpha,m}(x)(u(\sigma_\alpha x) - u(x)).$$

By Theorem 3.3 in ([5], p. 379) for any $\nu \in \mathcal{P}(E)$ the $\mathbf{D}_{\mathbb{R}^n}[0,\infty)$ martingale problem $(\mathcal{A}_m, \nu)$ is well-posed. By Theorem 6.1 in ([5], p. 216) the stopped $\mathbf{D}_{\mathbb{R}^n}[0,\infty)$ martingale problem $(\mathcal{A}_m, \nu, B_m)$ is well-posed. Since for any $u \in \mathbf{C}_K^\infty(\mathbb{R}^n)$, $x \in B_{m+1}$ $\mathcal{A}_m u(x) = \mathcal{L}_{k,l} u(x)$ and using (29), $X_{\cdot \wedge T_m}$ is a solution of the stopped $\mathbf{D}_{\mathbb{R}^n}[0,\infty)$ martingale problem $(\mathcal{A}_m, \nu, B_m)$. Hence, the distribution of $(X_{t \wedge T_m})_{t \geq 0}$ is uniquely determined. Using (28) one obtains that the distribution of $(X_t)_{t \geq 0}$ is uniquely determined. $\square$

Suppose that $R = \{\pm\alpha_1, \pm\alpha_2, \ldots, \pm\alpha_m\}$, $R_+ = \{\alpha_1, \alpha_2, \ldots, \alpha_m\}$. For any $i = 1, \ldots, m$ define $R^i = \{\pm\alpha_1, \pm\alpha_2, \ldots, \pm\alpha_i\}$, $R_+^i = \{\alpha_1, \alpha_2, \ldots, \alpha_i\}$, $R^0 = R_+^0 = \emptyset$.

**Theorem 19.** *Let $X$ be the Dunkl process, with extended generator given by (2), such that $X_0 = x \in E$ a.s. and $X^W$ – its radial part. Suppose that for any $\alpha \in R_+$, $k(\alpha) \geq \frac{1}{2}$.*

(i) *For $i = 1, \ldots, m$ there exist Poisson processes $N^i$ with intensity $k(\alpha_i)$ respectively and processes $Y^i$, defined recursively by*

$$Y_t^0 := X_t^W$$

*and*

$$Y_{\tau_\cdot^i}^i := Y_{\tilde{\tau}_\cdot^{i-1}}^{i-1} *_{\alpha_i} N^i, \tag{30}$$

*where*

$$A_t^i := \int_0^t \frac{ds}{(Y_s^i \cdot \alpha_i)^2}, \qquad \tau_t^i := \inf\{s \geq 0 \,|\, A_s^i > t\} \tag{31}$$

*and*

$$\tilde{A}_t^{i-1} := \int_0^t \frac{\mathrm{d}s}{(Y_s^{i-1} \cdot \alpha_i)^2}, \qquad \tilde{\tau}_t^{i-1} := \inf\{s \geq 0 \mid \tilde{A}_s^{i-1} > t\}, \tag{32}$$

such that for any $t > 0$ $A_t^i < +\infty$ a.s., $\tau_t^i < +\infty$ a.s., $\tilde{A}_t^{i-1} < +\infty$ a.s., $\tilde{\tau}_t^{i-1} < +\infty$ a.s., and $(X_t)_{t \geq 0} \stackrel{(d)}{=} (Y_t^m)_{t \geq 0}$, i.e. $Y^m$ is a Dunkl process with extended generator given by (2).

(ii) *For any $i = 0, \ldots, m$ $Y^i$ is a Markov process with extended generator $\mathcal{G}^i$, such that for any $u \in \mathbf{C}_K^\infty(E)$*

$$\mathcal{G}^i u(x) = \frac{1}{2}\Delta u(x) + \sum_{\alpha \in R_+} k(\alpha) \frac{\nabla u(x) \cdot \alpha}{x \cdot \alpha} + \sum_{\alpha \in R_+^i} k(\alpha) \frac{u(\sigma_\alpha x) - u(x)}{(x \cdot \alpha)^2}. \tag{33}$$

(iii) *If for some $i = 1, \ldots, m$*

$$\sigma_{\alpha_i}(R^{i-1}) = R^{i-1}, \tag{34}$$

*then $Y^i$ can be given by*

$$Y_t^i = \sigma_{\alpha_i}^{N^i \int_0^t (Y_s^{i-1} \cdot \alpha_i)^{-2} \, \mathrm{d}s} Y_t^{i-1}, \tag{35}$$

*where $N^i$ is a Poisson process with intensity $k(\alpha_i)$ independent from $Y^{i-1}$.*

**Remark 20.** *Let $k(\alpha) \geq \frac{1}{2}$ for any $\alpha \in R_+$. Fix a Weyl chamber $C$ and $x \in C$ and suppose that $X_0^W = x$ a.s. From the previous section we now that $X_t^W \in C$ for any $t \geq 0$ a.s. Let $Y^i$, $i = 0, \ldots, m$ be defined by Theorem 19. Then $Y^i \in C_i$ for any $t \geq 0$ a.s., where $C_i$ are defined recursively as follows:*

$$C_0 = C,$$
$$C_{i+1} = C_i \cup \sigma_{\alpha_i}(C_i) \tag{36}$$

*and $C_m = \mathbb{R}^n$.*

**Remark 21.** *The decomposition (30) depends on the way one enumerates the elements of R. Different enumerations lead to different skew-product decompositions of the Dunkl process.*

**Remark 22.** *One has $\sigma_{\alpha_m}(R^{m-1}) = R^{m-1}$ and $\sigma_{\alpha_1}(R^0) = R^0$. Therefore $Y^1$ and $Y^m$ can always be taken in the form (35).*

**Remark 23.** *There is some analogy between this decomposition and the skew-product decomposition of Brownian motion on the n-dimensional sphere in terms of the Legendre processes and Brownian motion on $(n-1)$-dimensional sphere in ([11], 7.15): one can iterate the skew-product decomposition of Brownian motion on n-dimensional sphere in order to get Brownian motion on $(n-2)$-dimensional sphere, then Brownian motion on $(n-3)$-dimensional sphere, etc.*

**Example 24.** *Take $R = B_2$ and $R_+ = \{\alpha_1, \ldots, \alpha_4\}$, where $\alpha_1 := e_1 - e_2$, $\alpha_2 := e_1 + e_2$, $\alpha_3 := e_1$, $\alpha_4 := e_2$, and $C = \{x = (x_1, x_2) \in \mathbb{R}^2 \mid x_2 > 0 \text{ and } x_1 > x_2\}$. Then the condition (34) is true for $i = 1, \ldots, 4$ and we obtain the skew-product decomposition (35).*

**Example 25.** *The condition (34) is not always satisfied. As a counterexample one can take $R = A_2$.*

**Proof of Theorem 19.** The proof is done by induction. Let us fix a Weyl chamber $C$ of $W$ and consider $X^W$ – the radial part of the Dunkl process. For any $x \in E$ there is one and only one $w_x \in W$ such that $w_x(x) \in C$. Define $\tilde{X}^W$ started at $x \in E$ by $w_x(X^W)$, where $X^W$ is the radial part of Dunkl process started at $w_x(x)$. One can see that the extended radial part $\tilde{X}^W$ is a Feller process on $E$ which is a.s. continuous with extended generator given by,

$$\mathcal{L}_k^W u(x) = \frac{1}{2}\Delta u(x) + \sum_{\alpha \in R_+} k(\alpha)\frac{\nabla u(x) \cdot \alpha}{x \cdot \alpha}, \qquad (37)$$

($u \in \mathbf{C}_K^\infty(E)$). Note also that $\pi(\tilde{X}^W)$ is a.s. continuous. Let $P_x$ be the law of $X^W$ on $\mathbf{C}_E[0,\infty)$. For any $\nu \in \mathcal{P}(E)$, let $P_\nu$ be the law on $\mathbf{C}_E[0,\infty)$ given by

$$P_\nu(\cdot) = \int_E \nu(\mathrm{d}x) P_x(\cdot).$$

Then the process $\bar{X}^W$ with the law $P_\nu$ is a solution of the $\mathbf{D}_E^\pi[0,\infty)$ martingale problem $(\mathcal{G}^0, \nu)$ (in order to recover the radial Dunkl process $X^W$ started at $x \in C$, one poses $\nu = \delta_x$). By Lemma 18 this solution is unique. Furthermore $\bar{X}^W$ is a Markov process.

For any $\nu \in \mathcal{P}(E)$ let $Y^{j-1}$ be the unique solution of the $\mathbf{D}_E^\pi[0,\infty)$ martingale problem $(\mathcal{G}^{j-1}, \nu)$. Since the paths of $Y^{j-1}$ are in $\mathbf{D}_E^\pi[0,\infty)$ the process $\log(\alpha_j \cdot Y^{j-1})^2$ is càdlàg in $\mathbb{R}$. Therefore, for any $t > 0$ $\log(\alpha_j \cdot Y^{j-1})^2$ is uniformly bounded on $[0,t]$ and $(\alpha_j \cdot Y^{j-1})^2$ is uniformly bounded from zero on $[0,t]$. Hence

$$\tilde{A}_t^{j-1} = \int_0^t \frac{\mathrm{d}s}{(Y_s^{j-1} \cdot \alpha_j)^2} < +\infty \quad \text{a.s.}$$

Next

$$\tilde{A}_t^{j-1} = \int_0^t \frac{\mathrm{d}s}{(Y_s^{j-1} \cdot \alpha_j)^2} \geq \frac{1}{2} \int_0^t \frac{\mathrm{d}s}{\|Y_s^{j-1}\|^2}.$$

But from (33) for any $f \in \mathbf{C}_K^\infty((0,+\infty))$ denoting $\eta_t := \|Y_t^{j-1}\|^2$

$$f(\eta_t) - f(\eta_0) - \int_0^t \mathrm{d}s[2\eta_s f''(\eta_s) + (n+2\gamma)f'(\eta_s)]$$

is a martingale. Since $\eta$ is the square of an $(n+2\gamma)$-dimensional Bessel process (started at $\eta_0$ with the law $\nu$)

$$\int_0^t \frac{1}{\eta_s}\mathrm{d}s \to \infty, \quad \text{as } t \to \infty \quad \text{a.s.}$$

Hence $\tilde{A}_t^{j-1} \to \infty$, as $t \to \infty$, and $\tilde{\tau}_\cdot^{j-1} < +\infty$ a.s. Finally $\tilde{\tau}_\cdot^{j-1}$ is a.s. strictly increasing continuous time-change with inverse $\tilde{A}_\cdot^{j-1}$. As in ([22], p. 168) for any $u \in \mathbf{C}_K^\infty(E)$

$$M_t^u := u(Y_{\tilde{\tau}_t^{j-1}}^{j-1}) - u(Y_0^{j-1}) - \int_0^t (Y_{\tilde{\tau}_s^{j-1}}^{j-1} \cdot \alpha_j)^2 \mathcal{G}^{j-1} u(Y_{\tilde{\tau}_s^{j-1}}^{j-1}) \, \mathrm{d}s \qquad (38)$$

is a local martingale. Since $(M_t^u)_{t\geq 0}$ is a.s. bounded on bounded time intervals it is a martingale. Hence, $Y_{\tilde{\tau}_\cdot^{j-1}}^{j-1}$ is a solution of the $\mathbf{D}_E^\pi[0,\infty)$ martingale problem $((\cdot \cdot \alpha_j)^2 \mathcal{G}^{j-1}, \nu)$. Next by Lemma 17 there exists a Poisson process $N^j$ with intensity $k(\alpha_j)$ such that the process $Z^j = Y_{\tilde{\tau}_\cdot^{j-1}}^{j-1} *_{\alpha_j} N^j$ is a solution of the $\mathbf{D}_E^\pi[0,\infty)$ martingale problem $(\mathcal{A}^j, \nu)$, where for any $u \in \mathbf{C}_K^\infty(E)$

$$\mathcal{A}^j u(x) = (x \cdot \alpha_j)^2 \mathcal{G}^{j-1} u(x) + k(\alpha_j)(u(\sigma_{\alpha_j} x) - u(x)).$$

Note that $(\int_0^t (\alpha_j \cdot Z_s^j)^2 \, ds)_{t \geq 0}$ is an a.s. strictly increasing, finite, continuous time-change. Let $\bar{\tau} := \lim_{t \to +\infty} \int_0^t (\alpha_j \cdot Z_s^j)^2 \, ds$. For any $t < \bar{\tau}$, define $\tau(t)$ as a unique solution of

$$t = \int_0^{\tau(t)} (\alpha_j \cdot Z_s^j)^2 \, ds.$$

Then

$$\frac{d}{dt}\tau(t) = \frac{1}{(\alpha_j \cdot Z_{\tau(t)}^j)^2}.$$

Define now $Y_t^j := Z_{\tau(t)}^j$, for $t < \bar{\tau}$. Then

$$\tau(t) = \int_0^t \frac{ds}{(\alpha_j \cdot Y_s^j)^2},$$

$\tau(t) < +\infty$ for any $t < \bar{\tau}$ and $\tau(t) \to +\infty$, as $t \to \bar{\tau}$. Hence, there exists a sequence $\{\tau_n\} \subset (0, \bar{\tau})$, such that $\tau_n \to \bar{\tau}$ and $(\alpha_j \cdot Y_{\tau_n}^j)^2 \to 0$, as $n \to +\infty$. Since for any $u \in \mathbf{C}_K^\infty(E)$

$$u(Z_t^j) - u(Z_0^j) - \int_0^t ((Z_s^j \cdot \alpha_j)^2 \mathcal{G}^{j-1} u(Z_s^j) + k(\alpha_j)(u(\sigma_{\alpha_j} Z_s^j) - u(Z_s^j))) \, ds$$

is a martingale, as for (38),

$$u(Y_t^j) - u(Y_0^j) - \int_0^t (\mathcal{G}^{j-1} u(Y_s^j) + k(\alpha_j)(Y_s^j \cdot \alpha_j)^{-2}(u(\sigma_{\alpha_j} Y_s^j) - u(Y_s^j))) \, ds, \quad t < \bar{\tau},$$

is a martingale. Denote $Y^W := \pi(Y^j)$, then for any $u \in \mathbf{C}_K^\infty(C)$

$$u(Y_t^W) - u(Y_0^W) - \int_0^t \mathcal{L}_k^W u(Y_s^W) \, ds, \quad t < \bar{\tau},$$

is a martingale. This shows that $Y^W$ is a radial Dunkl process up to time $\bar{\tau}$. Since $(\alpha_j \cdot Y_{\tau_n}^j)^2 \to 0$, as $n \to +\infty$, $\lim_{n \to +\infty} Y_{\tau_n}^W \notin C$ and for any $\alpha \in R_+$ $k(\alpha) \geq \frac{1}{2}$, by the results of the previous section $\tau_n \to \infty$ a.s. and $\bar{\tau} = +\infty$ a.s. Hence

$$\int_0^t \frac{ds}{(\alpha_j \cdot Y_s^j)^2} \to +\infty, \quad \text{as } t \to +\infty$$

and for any $t > 0$ $A_t^j < +\infty$, $\tau_t^j < +\infty$ a.s. Furthermore $Y^j$ is a solution of $\mathbf{D}_E^{\bar{\tau}}[0, \infty)$ martingale problem $(\mathcal{G}^j, \nu)$. Now by Lemma 18 the $\mathbf{D}_E^{\bar{\tau}}[0, \infty)$ martingale problem $(\mathcal{G}^j, \nu)$ is well-posed. Since the $\mathbf{D}_E^{\bar{\tau}}[0, \infty)$ martingale problem $(\mathcal{G}^j, \nu)$ is well-posed for any $\nu \in \mathcal{P}(E)$, by the same argument as in Theorem 4.2(a) in ([5], p. 184), $Y^j$ is a Markov process.

Finally when $j = m$ one obtains that $Y^m$ is the unique solution of the $\mathbf{D}_E^{\bar{\tau}}[0, \infty)$ martingale problem $(\mathcal{G}^m, \nu)$, but for any $u \in \mathbf{C}_K^\infty(E)$

$$\mathcal{G}^m u(x) = \mathcal{L}_k u(x).$$

Therefore $Y^m$ is a Dunkl process.

In order to get (iii) note that from (34) for any $u \in \mathbf{C}_K^\infty(E)$ and $x \in E$

$$\mathcal{G}^{i-1}(u \circ \sigma_{\alpha_i})(x) = \mathcal{G}^{i-1} u(\sigma_{\alpha_i}(x)),$$

and by Lemma 17 one can take $Y^i_{\tau^i_{\cdot}} = Z^i_{\cdot}$, where $\tau^i_t = \int_0^t (\alpha_i \cdot Z^i_s)^2 \, ds$ and

$$Z^i_{\cdot} = \sigma_{\alpha_i}^{N^i_{\cdot}} Y^{i-1}_{\tilde{\tau}^{i-1}_{\cdot}},$$

where $N^i$ is a Poisson process with intensity $k(\alpha_i)$ independent from $Y^{i-1}_{\tilde{\tau}^i}$ (hence, independent from $Y^{i-1}$). But

$$(\sigma_{\alpha_i} x \cdot \alpha_i)^2 = ((\alpha_i \cdot x) - (\alpha_i \cdot x)(\alpha_i \cdot \alpha_i))^2 = (x \cdot \alpha_i)^2$$

and

$$\tau^i_t = \int_0^t (\alpha_i \cdot Z^i_s)^2 \, ds = \int_0^t (\alpha_i \cdot Y^{i-1}_{\tilde{\tau}^{i-1}_s})^2 \, ds. \tag{39}$$

Finally differentiating the equality

$$t = \tilde{A}^{i-1}_{\tilde{\tau}^{i-1}_t}$$

one gets that

$$\frac{d}{dt} \tilde{\tau}^{i-1}_t = (\alpha_i \cdot Y^{i-1}_{\tilde{\tau}^{i-1}_t})^2$$

and from (39) $\tau^i_{\cdot} = \tilde{\tau}^{i-1}_{\cdot}$. Hence,

$$Y^i_t = \sigma_{\alpha_i}^{N^i_{\tilde{A}^{i-1}_t}} Y^{i-1}_t. \qquad \square$$

From Theorem 19 if (34) and (40) hold for some $j$ we can deduce the following relationship between the semigroups of $Y^j$ and $Y^{j-1}$.

**Proposition 26.** *Under the conditions of Theorem 19 suppose that for a certain $j$ (34) holds. Fix a Weyl chamber $C$. Let $Y^0_0 \in C$ a.s., then $Y^{j-1} \in C_{j-1}$ and $Y^j \in C_j$ a.s., where $(C_i)$ are defined by (36). Suppose that $P^{j-1}_t(x, dy)$ and $P^j_t(x, dy)$ are the semi-groups of $Y^{j-1}$ and $Y^j$ respectively, and*

$$C_{j-1} \cap \sigma_{\alpha_j}(C_{j-1}) = \emptyset, \tag{40}$$

*then for any $x, y \in C_{j-1}$*

$$P^{j-1}_t(x, dy) = P^j_t(x, dy) + P^j_t(x, \sigma_{\alpha_j}(dy)). \tag{41}$$

**Proof.** By Remark 20

$$C_j = C_{j-1} \cup \sigma_{\alpha_j}(C_{j-1}).$$

For any bounded measurable $f : C_{j-1} \to \mathbb{R}$ define $g : C_j \to \mathbb{R}$ such that for any $x \in C_j$

$$g(x) := f(x) \mathbb{I}_{\{x \in C_{j-1}\}} + f(\sigma_{\alpha_j}(x)) \mathbb{I}_{\{x \in \sigma_{\alpha_j}(C_{j-1})\}},$$

then for any $x \in C_{j-1}$

$$\mathbb{E}_x g(Y^j_t) = \int_{C_{j-1}} f(y) P^j_t(x, dy) + \int_{\sigma_{\alpha_j}(C_{j-1})} f(\sigma_{\alpha_j}(y)) P^j_t(x, dy)$$

$$= \int_{C_{j-1}} f(y) (P^j_t(x, dy) + P^j_t(x, \sigma_{\alpha_j}(dy))). \tag{42}$$

On the other hand

$$\mathbb{E}_x g(Y_t^j) = \mathbb{E}_x f(Y_t^{j-1}) = \int_{C_{j-1}} f(y) P_t^{j-1}(x, \mathrm{d}y). \tag{43}$$

Comparing (42) and (43) one obtains (41). □

Using Lemma 18 and the analog of Theorem 19 (one should change $k$ to $k'$ in the proper places) one can introduce a slightly more general class of Markov processes than Dunkl processes – the $(k, k')$-Dunkl processes. These processes are introduced in dimension 1 in [10] and in general dimensions in [7] (see also [15]) and are two parameter analogs of Dunkl processes. They are characterized by their extended generator, for any $u \in \mathbf{C}_K^\infty(E)$,

$$\mathcal{L}_{k,k'} u(x) = \frac{1}{2}\Delta u(x) + \sum_{\alpha \in R_+} k(\alpha) \frac{\nabla u(x) \cdot \alpha}{x \cdot \alpha} - \sum_{\alpha \in R_+} k'(\alpha) \frac{u(x) - u(\sigma_\alpha x)}{(x \cdot \alpha)^2}, \tag{44}$$

where $k : R \to [\frac{1}{2}, \infty)$ and $k' : R \to [0, \infty)$ are two multiplicity functions invariant by the finite reflection group $W$ associated with $R$. The rest of the notations is the same as in (2). Denote such processes by $X^{(k,k')}$. It is simple to see that the radial part of $X^W = \pi(X^{(k,k')})$ is the same as for the Dunkl process $X^{(k)}$. Note that by Lemma 18 $X^{(k,k')}$, started at $x \in E$, is a unique solution of the $\mathbf{D}_E^\pi[0, \infty)$ martingale problem $(\mathcal{L}_{k,k'}, \delta_x)$.

## Acknowledgments

The author thanks Marc Yor whose ideas and help were crucial for this work. The author thanks an anonymous referee for valuable comments and suggestions.